\documentclass[11pt]{article}
\usepackage[latin1]{inputenc}   

\usepackage{amsfonts,amsmath,layout}

\newcommand{\N}{\mathbb N}

\newcommand{\R}{\mathbb R}

\def\QED{\hfill $\; \Box$\medskip}

\def\P{\mathbb P}

\def\A{\mathcal A}

\newcommand{\be}{\begin{equation}}
\newcommand{\ee}{\end{equation}}
\def\1{{\bf 1}}

\def\V{{\bf V}}

\newcommand{\B}{\mathcal{B}}

\def\ep{\epsilon}

\def\Dt0{{\bf D}(t_0)}

\def\E{{\bf E}}

\def\P{{\bf P}}
\def\bV{{\bf V}}

\def\Int{\displaystyle\int}

\def\to{\rightarrow}

\def\ds{\displaystyle}

\newcommand{\Vt}{\mathcal{V}}
\newcommand{\U}{\mathcal{U}}
\def\ind{{\bf 1}\hskip-2.5pt{\rm l}}

\newcommand{\ba}{\[\begin{array}{rl}}
\newcommand{\ea}{\end{array}\]}
\newcommand{\bea}{\begin{eqnarray}}
\newcommand{\eea}{\end{eqnarray}}
\newcommand{\beaa}{\begin{eqnarray*}}
\newcommand{\eeaa}{\end{eqnarray*}}
\newtheorem{Theorem}{Theorem}[section]
\newtheorem{Definition}[Theorem]{Definition}
\newtheorem{Proposition}[Theorem]{Proposition}
\newtheorem{Lemma}[Theorem]{Lemma}
\newtheorem{Corollary}[Theorem]{Corollary}
\newtheorem{Remark}[Theorem]{Remark}
\newtheorem{Remarks}[Theorem]{Remarks}

\newcommand{\AR}{{\cal A}}
\newcommand{\BR}{{\cal B}}

\newcommand{\FR}{{\cal F}}

\newcommand{\JR}{{\cal J}}

\newcommand{\MR}{{\cal M}}
\newcommand{\NR}{{\cal N}}

\newcommand{\UR}{{\cal U}}
\newcommand{\VR}{{\cal V}}

\begin{document}

\title{Games with incomplete information in continuous time and 
for continuous types.}

\author{Pierre Cardaliaguet$^{\tiny (1)}$ and Catherine Rainer$^{\tiny (2)}$\\
$\;$\\
\small {\tiny (1)} Universit\'e Paris Dauphine,\\
\small {\tiny (2)} Universit\'e de Bretagne Occidentale, 6, avenue Victor-le-Gorgeu, B.P. 809, 29285 Brest cedex, France\\
 \small e-mail : 
cardaliaguet@ceremade.dauphine.fr, Catherine.Rainer@univ-brest.fr}
  \maketitle

\noindent {Abstract : We consider a two-player zero-sum game with integral payoff and with incomplete information on one side, where the payoff is chosen among a continuous set of possible payoffs. We prove that the value function of this game is solution of an auxiliary optimization problem over a set of measure-valued processes. Then we use this equivalent formulation to characterize the value function as the viscosity solution of a special type of a Hamilton-Jacobi equation. This paper generalizes the results of a previous work  of the authors \cite{cr2}, where only a finite number of possible payoffs is considered. 

}
\vspace{3mm}

\vspace{3mm}

\noindent{\bf A.M.S. classification :} 91A05, 91A23, 60G60
\vspace{3mm}

\section{Introduction.}

In this paper, we investigate a continuous time two-player zerosum game with lack of information on one side:   according to a given probability, some payoff is chosen randomly among a family of payoffs depending on a parameter taking its values in $\R^N$, and is communicated to one player (say Player 1), Player 2 being only aware of the underlying  probability. Then Player 1 want to minimize, Player 2 to maximize the payoff.  \\
Games with incomplete informations (where some players have private informations) were introduced in the framework of repeated games by Aumann and Maschler in the nineteen sixties (see \cite{AumannMaschler}) and are, since then, the aim of many works, see for instance \cite{demeyerrozenberg}, \cite{mertenszamir}, \cite{sorin} and the references therein.\\
The first author of the present paper introduced the continuous time case in \cite{c1}. It is proven there that, under suitable assumptions, differential games with asymmetric informations have a value which can be characterized as a solution in some dual, viscosity sense of a Hamilton-Jacobi-Isaacs equation. The result was generalized to stochastic differential games in \cite{cr1}. \\
A crucial point in games with incomplete information is that, the players observing each other all along the game, the non informed player will try to guess his missing information through the actions of his opponent. This implies that the informed player, instead of fully using his information, needs to reveal it sparingly and to hide its revelation by using random strategies to this aim. In the framework of repeated games, the behavior of the informed player is therefore strongly linked with a martingale with values in the set of probability measures, called ``martingale of beliefs'', which describes the belief the non informed player has on the chosen index. Rewriting the value of the game as an optimization problem over the set of martingales has been extensively used in the recent literature \cite{DeM, Gen}. When the set of possible payoffs is finite (say: $I$ possible payoffs), this set of probabilities over $\{1, \dots, I\}$ can be assimilated to the simplex $\Delta(I)$, subset of $\R^I$. In the case of continuous times games without dynamic, and in analogy to the repeated games, we proved in \cite{cr2}, that the value function of the game is solution of an optimization problem over this set of $\R^I$-valued martingales. This result is generalized to stochastic differential games by Gr\"un \cite{gruen}.\\
In the present paper we investigate the case of continuous games without dynamics, with a continuum of possible cost functions. A crucial aspect of games with a finite number of possible payoffs, is that the set of probability-measures can be assimilated to $\R^I$. Therefore the value function depend on two parameters $t\in[0,T]$, the starting time and $p\in\Delta(I)$ the initial probability and, even if the place taken by the parameter $p$  in the Hamilton-Jacobi-Isaacs equation is not standard, it is still an equation in $\R\times\R^I$. The equivalent formulation of the value function as a control problem over a set of martingales In the case of a continuum of possible payoffs. 

The main results of this paper are, on the one hand, a representation formula of the value in terms of a minimization problem over a set of martingale measures and the existence of an optimal martingale measure. On another hand, we also characterize the value as the (dual) solution of a Hamilton-Jacobi equation. 

\section{Notations and preliminaries.}

Let us first introduce the model we will study throughout the paper. We consider a continuous time game starting at some initial time $t_0>0$ and ending at a terminal time $T>t_0$. The players try to optimize an integral payoff of the form
$$
\int_{t_0}^T \ell (x, s, u(s), v(s)) \ ds 
$$
where $(u(s))_{s\in [t_0,T]}$ is the control played by the first player while  $(v(s))_{s\in [t_0,T]}$ is the control played by the second player. The first player is minimizing, the second one maximizing. The main point in the above payoff function is its dependence with respect to a parameter $x$. The game is played as follows: at time $t_0$ the parameter $x$ is chosen by nature according to some probability $\bar m$ and  the result is announced to the Player 1, but not to the Player 2. Then the players choose their respective controls in order to optimize the payoff. The map $\ell$ as well as the measure $\bar m$ are public knowledge. 

To fix the ideas, we will assume throughout the paper that the parameter $x$ belongs to some Euclidean space $\R^N$. Although this assumption is not really important for the existence of a value, the representation Theorem  \ref{altform} or the characterization of the value given in Proposition \ref{largest}, it will be crucial for the existence of an optimal martingale measure (Theorem \ref{theo:existsoptimal}) or a sharper characterization (Proposition \ref{Prop:Charac2}).   

Moreover, to get some compactness ensuring the existence of an optimal martingale measure, we will have to assume that the measure $\bar m$ has some finite moment. Again, to fix the ideas, we will assume the finiteness of the second moment (however, the reader will easily notice that any moment larger than 1 would do the job).
Let $\Delta(\R^N)$ be the set of Borel probability measures on $\R^N$ and  $\Pi_2$ be the subset of $\Delta(\R^N)$ with a finite second order moment:
$$
\Int_{\R^N}|x|^2dm(x) <+\infty \qquad \forall m\in \Pi_2\;.
$$
It will be convenient to endow $\Pi_2$ with the Monge-Kantorovich distance ${\bf d}_1$:
$$
{\bf d}_1(m,m')= \inf_{\gamma\in \Pi(m,m')}  \Int_{\R^{2N}} |x-y|d\gamma(x,y)
$$
where $\Pi(m,m')$ is the set of Borel probability measures on $\R^{2N}$ such that $\gamma(A\times \R^N)=m(A)$
and $\gamma(\R^N\times A)=m'(A)$ for any Borel set $A\subset \R^N$. Recall that the above minimization problem has at least one minimum
and we denote by $\Pi_{\rm opt}(m,m')$ the set of such minima for ${\bf d}_1$. Again there is some arbitrariness in the choice of this distance, but it is useful here to fix the notations. 

Let $T>0$ be a finite time horizon, and $U$ and $V$ two compact subsets of some finite dimensional space.  For all $0\leq t_0\leq t_1\leq T$ we denote by 
$\UR(t_0,t_1)=\{ u:[t_0,t_1]\rightarrow U\  \mbox{measurable}\}$
  the set of controls with values in $U$. 
We endow $\U(t_0,t_1)$ with the $L^1$ distance
$$
d_{\U(t_0,t_1)}(u_1,u_2)= \Int_{t_0}^{t_1} |u_1(s)-u_2(s)|ds \qquad \forall u_1,u_2\in \U(t_0,t_1)\;,
$$
and with the Borel $\sigma-$algebra associated with this distance. Recall that $\U(t_0,t_1)$ is then a Polish space (i.e., a complete separable
metric space). The set $\Vt(t_0,t_1)$ of Lebesgue
measurable maps $v:[t_0,t_1]\to V$ is defined in a symmetric way and 
endowed with the $L^1$ distance and with the associate Borel $\sigma-$algebra. 
We write $\UR(t_0)$ and $\Vt(t_0)$ for  $\UR(t_0,T)$ and $\Vt(t_0,T)$.

The instantaneous reward \[ \ell:\R^N\times[0,T]\times U\times V\rightarrow \R\]
is assumed to be continuous in all variables, Lipschitz in $(x,t)$ uniformly in $u$ and $v$, and bounded.\\
Fix now $(t_0,m)\in[0,T]\times\Pi_2$. For all $(u,v)\in\UR(t_0)\times\VR(t_0)$, we define the payoff
\[ \JR(t_0,m,u,v)=\int_{\R^N}\int_{t_0}^T\ell(x,s,u(s),v(s))\ ds\ m(dx).\]
Note that $J(t_0,m,u,v)$ is the mean with respect to the probability $m$ of a family of running costs depending on the parameter $x\in\R^N$. 
The control $u$ represents the action of Player 1, $v$ the action of Player 2. In this game, Player 1 will try to minimize, Player 2 to maximize the payoff.\\
In order to give a probabilistic interpretation of the game, we now introduce a probability space $(\Omega_0,\FR_0,P_0)$. If $X$ is a random variable of distribution $m$, we can rewrite the payoff 
\[ \JR(t_0,m,u,v)=E_0\left[\int_{t_0}^T\ell\left(X,s,u(s),v(s)\right)ds\right].\]
The next step is to define the strategies of the Players. We have to model the facts that
\begin{itemize}
\item Player 1 knows the index $x$ before the game starts, but not Player 2,
\item each Player observe the action of his opponent,
\item the Players can act randomly.
\end{itemize}

\noindent The strategies of Player 1 	are defined as follows:
For some fixed subdivision $\Delta:t_0<t_1<\ldots<t_n=T$, let $\A_\Delta(t_0)$ be the set of Borel measurable maps $\alpha:\Vt(t_0)\to \U(t_0)$ which are
nonanticipative with delay with respect to the subdivision $\Delta$ :  
for any $v_1,v_2\in\Vt(t_0)$, if $v_1\equiv v_2$ a.e. on $[t_0,t_i]$ for some $i\in \{0, \dots, n-1\}$, then
$\alpha(v_1)\equiv \alpha(v_2) \;  {\rm on }\; [t_0,t_{i+1}]
$.  
We set $\AR(t_0)=\cup_{\Delta}\AR_\Delta(t_0)$ and endow $\A(t_0)$ with the distance
$$
d(\alpha_1, \alpha_2)= \sup_{v\in \Vt(t_0)} d_{\U(t_0)}(\alpha_1(v),\alpha_2(v))
$$

\noindent An element of  $\BR(t_0)$, the set of strategies for Player 2, is a nonanticipative, Borel measurable map $\beta:\U(t_0)\to \Vt(t_0)$ such that, for all $t\in[t_0,T]$,  for any $u_1,u_2\in\U(t_0)$, there exists $\delta>0$ such that, if $u_1\equiv u_2$ a.e. on $[t_0,t]$, then
$\beta(u_1)\equiv \beta(u_2) \;  {\rm on }\; [t_0,t+\delta]$. (i.e, in contrast to the strategies of Player I, Player II is allowed to adapt his delay belonging what he is answering to).
\\

Due to the delays the Players have to respect, we get, as in \cite{cr1, cr2}:
\begin{Lemma} 
\label{fixpoint}
There exists a Borel measurable  map $\Phi: \A(t_0)\times \B(t_0)\to\U(t_0)\times \Vt(t_0)$ such that, for all pair $(\alpha,\beta)\in\A(t_0)\times \B(t_0)$ and  $(u,v)\in\U(t_0)\times \Vt(t_0)$, 
\begin{equation}
\label{fix}
(u,v):=\Phi(\alpha,\beta) \qquad \Leftrightarrow \qquad\left[ \ \alpha(v)=u\; {\rm and }\; \beta(u)=v \; {\rm a.e.} \  \right]
\end{equation}
\end{Lemma}

\begin{Definition} Fix  some other  sufficiently large probability spaces $(\Omega_1,\FR_1,P_1)$ and $(\Omega_2,\FR_2,P_2)$. A strategy for Player I for the initial time $t_0\in[0,T]$  is a measurable application $\alpha$  from $\R^N\times\Omega_1$  to some $\A_\Delta(t_0)$.
We denote by $\tilde\A_X(t_0)$ the set of these strategies. \\
A strategy for Player II is a random variable on $(\Omega_2,\FR_2,P_2)$ with values in $\B(t_0)$.
We denote by $\tilde\BR(t_0)$ the set of strategies for Player II.
\end{Definition}

\noindent We set $(\Omega,\FR,P)=(\Omega_0\times\Omega_1\times\Omega_2,\FR_0\otimes\FR_1\otimes\FR_2,P_0\otimes P_1\otimes P_2)$ and denote by $E[\cdot]$ the corresponding expectation.\\

\noindent Fix $t_0\in[0,T]$. For any $(\alpha,\beta)\in\tilde\A_X(t_0)\times\tilde\B(t_0)$, we are now able to define
the payoff of the two strategies by
$$
\begin{array}{rl}
 \JR(t_0,\alpha,\beta,m)\; := & \ds \int_{\R^N\times \A(t_0)\times \B(t_0)}\Int_{t_0}^T\ell(x,t,\Phi(\alpha(x,\omega_1),\beta(\omega_2))(t))
dt\ dP_1(\omega_1)dP_2(\omega_2)dm(x)\\
= & \ds E\left[\Int_{t_0}^T\ell(X,t,\Phi(\alpha(X),\beta)(t))dt\right]
\end{array}
$$

The upper and lower value functions of the game are defined by
\[ \bV^+(t_0,m)=\inf_{\alpha\in\tilde\AR_X(t_0)}\sup_{\beta\in\tilde\BR(t_0)}\JR(t_0,\alpha,\beta,m),\]
and
\[ \bV^-(t_0,m)=\sup_{\beta\in\tilde\BR(t_0)}\inf_{\alpha\in\tilde\AR_X(t_0)}\JR(t_0,\alpha,\beta,m).\]
Note that we also have
\[ \bV^+(t_0,m)=\inf_{\alpha\in\tilde\AR_X(t_0)}\sup_{\beta\in\BR(t_0)}\JR(t_0,\alpha,\beta,m).\]

\begin{Lemma}
\label{lipschitzm}
The value functions $\bV^+$ and $\bV^-$ are  Lipschitz continuous. 
\end{Lemma}

\noindent{\bf Proof:} We only explain the proof for $\bV^+$, the arguments for $\bV^-$ being symmetrical. Let us first check the Lipschitz continuity in $m$.  Let $m_1,m_2\in \Pi_2$ and $\gamma\in \Pi_{\rm opt}(m_1,m_2)$. Let us desintegrate $\gamma$ with respect to $m_1$:
$d\gamma(x,y)=d\gamma_x(y)dm_1(x)$. \\
For $\epsilon>0$, let $\alpha_2\in\tilde\AR_X(t_0)$ be $\epsilon$-optimal for $\bV^+(t_0,m_2)$:
\[ \sup_{\beta\in\BR(t_0)}\JR(t_0,\alpha_2,\beta,m_2)\leq \bV^+(t_0,m_2)+\epsilon.\]
Let $\xi:\R^N\times\Omega_1\rightarrow\R^N$ be a measurable map such that, for all $x\in\R^N$, the random variables $(\xi(x,\cdot), x\in\R^N)$ are independent of $(\alpha_2(x),x\in\R^N)$ and the distribution of each $\xi(x,\cdot)$ is
$\gamma_x$.\\
Then the map $(\omega,x,v)\rightarrow\alpha_1(\omega,x,v):=\alpha_2(\omega,\xi(x,\omega),v)$ defines a strategy, i.e., an element of $\tilde\AR_X(t_0)$. We have, for all $\beta\in\BR(t_0)$,
\[\begin{array}{rl}
 \JR(t_0,\alpha_1,\beta,m_1)=&
\Int_{\R^N}E_1\left[\Int_{t_0}^T\ell(x,t,\Phi(\alpha_2(\xi(x)),\beta)(t))dt\right]dm_1(x)\\
\\
=&
\Int_{\R^N}dm_1(x)\int_{\R^N}d\gamma_x(y)E_1\left[\Int_{t_0}^T\ell(x,t,\Phi(\alpha_2(y),\beta)(t))dt\right]\\
\\
=&\Int_{\R^{2N}}d\gamma(x,y)E_1\left[\Int_{t_0}^T\ell(x,t,\Phi(\alpha_2(y),\beta)(t))dt\right]\\
\\
\leq&\Int_{\R^{2N}}d\gamma(x,y)E_1\left[\int_{t_0}^T \left(\ell(y,t,\Phi(\alpha_2(y),\beta)(t))+C|x-y|\right)dt\right]\\
\\
\leq &\JR(t_0,\alpha_2,\beta,m_2)+C{\bf d}_1(m_1,m_2)\\
\\
\leq & \bV^+(t_0,m_2)+\epsilon+C{\bf d}_1(m_1,m_2).
\end{array}\]
Therefore
\[ \bV^+(t_0,m_1)\leq \sup_{\beta\in\BR(t_0)}\JR(t_0,\alpha_1,\beta,m_1)\leq \bV^+(t_0,m_2)+\epsilon+C{\bf d}_1(m_1,m_2).\]
This proves the Lipschitz continuity of $\bV^+$ with respect to $m$, uniformly in $t$. \\

Next we prove the Lipschitz continuity in time. Fix $m\in\Pi_2$, $0\leq t_0\leq t_1\leq T$ and $\bar u\in U$.\\
For $\epsilon>0$, let $\alpha_1$ be $\epsilon$-optimal for $\bV^+(t_1,m)$.
We define a strategy in $\tilde\AR_X(t_0)$ by setting
\[ \alpha_0(x,v,\omega_1)(t)=
\left\{
\begin{array}{ll}
\bar u  &\mbox{ for } t\in[t_0,t_1)\\
\alpha_1(x,v_{|_{[t_1,T]}},\omega_1)&\mbox{ otherwise.}
\end{array}
\right.\]
For all $\beta_0\in\BR(t_0)$, we define $\beta_1\in\BR(t_1)$ by setting 
$\beta_1(u)= \beta_0((\bar u, u)) $, 
where $(\bar u,u)$ denotes the control equal to $\bar u$ on $[t_0,t_1)$ and to $u$ on $(t_1,T)$. 
Let $(u_1(x,\omega_1),v_1(x,\omega_1))=\Phi(\alpha_1(x,\omega_1),\beta_1)$ be the solution of the fixed point equation given in Lemma \ref{fixpoint}. 
Then, by definition of $\alpha_0$ and because $\beta_0$ is nonanticipative, $\Phi(\alpha_0(x,\omega_1),\beta_0)$ is equal to 
$(\bar u, \beta_0(\bar u))$ on $[t_0,t_1]$ and to $(u_1(x,\omega_1),v_1(x,\omega_1))$ on $[t_1,T]$.  
Hence, as the payoff is bounded by $C$,  \[\begin{array}{rl}
\JR(t_0,m,\alpha_0,\beta_0)=&
E\left[\int_{t_0}^T\ell(X,t,\Phi(\alpha_0(X),\beta_0)(t))dt\right]\\
\\
\leq & C(t_1-t_0)+\JR(t_1,m,\alpha_1,\beta_1)\\
\\
\leq &  C(t_1-t_0) + \bV^+(t_1,m)+\epsilon.
\end{array}\]
Therefore 
\[ \bV^+(t_0,m)\leq C(t_1-t_0) + \bV^+(t_1,m)\ .\]
The reverse inequality 
$$
\bV^+(t_0,m)\geq \bV^+(t_1,m)-C(t_1-t_0) 
$$
can be established in a symmetrical way. 
%
%
\QED

\begin{Lemma}
\label{convex}
The value functions and $\V^-$ and  $\bV^+$ are convex with respect to $m$.
\end{Lemma}

\noindent{\bf Proof:} The convexity of $\V^-$ is an easy consequence of the definition, because
$$
\bV^-(t_0,m)=\sup_{\beta\in\tilde\BR(t_0)}
\int_{\R^N} \inf_{\alpha\in\A(t_0)}\int_{\Omega_2} \int_0^T \ell (x, t, \Phi(\alpha,\beta(\omega_2)))dt \ dP_2(\omega_2)dm(x)\;,
$$
which is clearly convex with respect to $m$. 

Let us now check the convexity of $\V^+$. Fix $t_0\in[0,T]$ and let $m_1,m_2\in\Pi_2$ and $\lambda\in(0,1)$. Set $m=\lambda m_1+(1-\lambda) m_2$.\\
Since $m_1$ is absolutely continuous with respect to $m$, there exist $p_\lambda\in L^1(\R^N, \R_+; m)$ such that
\[ \lambda m_1(dx)=p_\lambda(x)m(dx) \]
with $p_\lambda\geq 0$ $m-$a.e. Remark that
\[ (1-\lambda)m_2(dx)=(1-p_\lambda(x))m(dx),\]
and therefore that $p_\lambda$ takes its values in $[0,1]$  $m-$a.e.
Now, for $\epsilon>0$, let $\alpha_1$ (resp. $\alpha_2$) $\in\tilde\AR_X(t_0)$ be $\epsilon$-optimal for $\bV^+(t_0,m_1)$ (resp. $\bV^+(t_0,m_2)$), and let $\alpha\in\tilde\AR_X(t_0)$ be such that, for all $x\in\R^N$, and all measurable, bounded $f:\A(t_0)\to \R$,
\[ E_1\left[f(\alpha(x))\right]=p_\lambda(x)E_1\left[f(\alpha_1(x))\right]+(1-p_\lambda(x))E_1\left[f(\alpha_2(x))\right].\]
For any $\beta\in\BR(t_0)$, it follows that
\[\begin{array}{l}
\JR(t_0,m,\alpha,\beta)=\Int_{\R^N}E_1\left[\int_{t_0}^T\ell(x,t,\Phi(\alpha(x),\beta)(t))dt\right]dm(x)\\
\\
=
\Int_{\R^N}\left( p_\lambda(x)E_1\left[\int_{t_0}^T\ell(x,t,\Phi(\alpha_1(x),\beta)(t))dt\right]
+(1-p_\lambda(x))E_1\left[\int_{t_0}^T\ell(x,t,\Phi(\alpha_2(x),\beta)(t))dt\right]\right)dm(x)\\
\\
=
\lambda \Int_{\R^N}E_1\left[\int_{t_0}^T\ell(x,t,\Phi(\alpha_1(x),\beta)(t))dt\right]dm_1(x)
+(1-\lambda)\Int_{\R^N}E_1\left[\int_{t_0}^T\ell(x,t,\Phi(\alpha_2(x),\beta)(t))dt\right]dm_2(x)\\
\\
\leq 
\lambda \bV^+(t_0,m_1)+(1-\lambda)\bV^+(t_0,m_2)+\epsilon.
\end{array}\]
Since this holds true for all $\beta\in\BR(t_0)$ and for any $\epsilon>0$, it follows that
\[ \bV^+(t_0,m)\leq \lambda \bV^+(t_0,m_1)+(1-\lambda)\bV^+(t_0,m_2).\]
\QED

%

\noindent We now introduce the Hamiltonian of the game :
\[ \forall (t,m)\in[0,T]\times\Pi_2,\; H(t,m)=\inf_{u\in U}\sup_{v\in V}\Int_{\R^N}\ell(x,t,u,v)dm(x),\]
and assume throughout the rest of the paper that Isaac's condition:
for all $t\in[0,T]$ and $m\in\Pi_2$,
\begin{equation}
 \label{isaacs}
\inf_{u\in U}\sup_{v\in V}\Int_{\R^N}\ell(x,t,u,v)dm(x)=\sup_{v\in V}\inf_{u\in U}\Int_{\R^N}\ell(x,t,u,v)dm(x).
 \end{equation}

\begin{Proposition}
\label{VpVm}
Under Isaac's condition (\ref{isaacs}), the game has a value: for all $(t_0,m)\in[0,T]\times\Pi_2$,
$\bV^+(t_0,m)=\bV^-(t_0,m)$.
\end{Proposition}

\noindent{\bf Proof:} The result has been proved for measures $m$ with finite support  and very slight different notion of strategies in \cite{cr2}. A careful examination of the proof shows that the result of \cite{cr2} also hold with the new definition of strategies. Let us explain how the condition on the measure can be removed. For an arbitrary measure $m\in\Pi_2$, let $(m_n)_{n\in\N}$ be a sequence of measures with finite support, such that ${\bf d}_1(m_n,m)\rightarrow_{n\rightarrow+\infty} 0$. Since, by lemma \ref{lipschitzm}, $\bV^+$ and $\bV^-$ are Lipschitz with respect to $m$, we have, for all $t_0\in[0,T]$, and all $n\in\N$,
\[ |\bV^+(t_0,m)-\bV^-(t_0,m)|\leq |\bV^+(t_0,m)-\bV^+(t_0,m_n)|+|\bV^-(t_0,m)-\bV^-(t_0,m_n)|\leq 2C{\bf d}_1(m,m_n),\]
where $C$ is a Lipschitz constant for $\bV^+$ and $\bV^-$.
The right-hand side term being  arbitrarily small, the result follows.
\QED

\noindent Let us complete this paragraph with a technical Lemma, which shows that the payoff associated to a control without delay may be approximated as close as needed by the same control but with delay.

\begin{Lemma}
\label{delaybeta}
Let $m\in\Pi_2$ and $\epsilon>0$.
For all $v\in\VR(0)$, there exists $\delta(v)>0$ such that, for all $u\in\UR(0)$,
\[ \left|\int_{\R^N}\Int_0^T\ell(x,s,u_s,v_s)dsdm(x)
- \int_{\R^N}\Int_0^T\ell(x,s,u_s,v_{(s-\delta(v))_+}dsdm(x)\right|\leq \epsilon.\] 
\end{Lemma}

\noindent{\bf Proof:}
Let $M$ be an upper bound for $|\ell(x,s,u,v)|$. 
Fix $\epsilon>0$ and let $B\in\BR(\R^N)$ a compact set such that $m(B)\geq 1-\epsilon/6MT$. 
Then we have, for all $\delta>0$,
\[ \left| \int_{\R^N}\int_0^T\ell(x,s,u_s,v_s)dsdm(x)-\int_{\R^N}\int_0^T\ell(x,s,u_s,v_{(s-\delta)_+})dsdm(x)\right|\leq \epsilon/3+I(\delta,\epsilon).\]
with
\[
 I(\delta,\epsilon):=\int_B\int_0^T\left|\ell(x,s,u_s,v_s)-\ell(x,s,u_s,v_{(s-\delta)_+})\right|dsdm(x)\]
Now, since $\ell$ is continuous,  it is uniformly continuous on the compact set  $B\times [0,T]\times U\times V$.
Therefore there exists $\eta>0$ such that, for all $(x,s,u)\in B \times [0,T]\times U$ and $v,v'\in V$ with $|v-v'|\leq \eta$, 
$|\ell(x,s,u,v)-\ell(x,s,u,v')|\leq\epsilon/3T$.\\
It follows that, for all $\delta>0$,

\[\begin{array}{rl}
 I(\delta,\epsilon)
\leq &  \epsilon/3+2M\lambda(\{s\in[0,T],| v_s-v_{(s-\delta)_+}|>\eta\}).
\\
\\
\leq &
 \epsilon/3+\frac{2M}\eta \int_0^T| v_s-v_{(s-\delta)_+}|ds,
\end{array}\]
where the last relation comes from the inequality of Bienaym\'e-Tchebychev.\\
Finally it is well known that, for any $v\in\VR(0)$,
\[ \lim_{\delta\searrow 0}\int_0^T| v_s-v_{(s-\delta)_+}|ds=0.\]
Therefore there exists $\delta_0>0$ such that, for all $\delta<\delta_0$, 
\[ \int_0^T| v_s-v_{(s-\delta)_+}|ds\leq \epsilon\eta/6M.\]
The result follows.
\QED

\section{An alternative formulation of the value function.}

In this section, we show that the upper value function can be represented in term of a minimization of a cost over a family of martingale measures. 
Many results on these processes used below can be found in Horowitz \cite{horowitz} or Walsh \cite{walsh} for instance.

For any $0\leq t_0\leq t_1\leq T$ and $m\in\Pi_2$, we introduce the set $\MR(t_0,t_1,m)$ of all measure-valued processes $(M_t)_{t_0\leq t\leq t_1}$ defined on a sufficiently large probability space $( \Omega, \FR, P)$ which satisfy  for all measurable, bounded functions $f:\R^N\rightarrow\R$,
\begin{itemize}
\item[(i)] $t\mapsto  \Int_{\R^N}f(x)dM_s(x)$ is   $ P$-a.s. c\`adl\`ag,
\item[(ii)] there exists a $\R^N$-valued random variable  $X$ of distribution $m$ such that, for all $t\in[t_0,t_1]$,  $ P$-a.s.,
\begin{equation}
\label{MX}
M_tf:= \Int_{\R^N}f(x)dM_t(x)=E\left[f(X)|\FR^M_t\right],
\end{equation}
where $(\FR^M_t)_{t\in[t_0,t_1]}$ is the right continuous, completed filtration generated by $(M_t(B), B\in\BR(\R^N))_{t\in[t_0,t_1]}$.
\end{itemize}

\begin{Remarks}
1) Let $M\in\MR(t_0,t_1)$.
\begin{itemize}
\item For all $t\in[t_0,t_1]$,  the random variable $M_t$ is $ P$-a.s. a probability measure  which satisfies
\[ E\left[M_t(B)\right]=m(B),\; B\in\BR(\R^N).\]
\item The process $(M_t)$ is a martingale measure:  for all measurable, bounded function $f$, $(M_tf)_{t\in[t_0,t_1)}$ is a martingale  in the filtration $(\FR^M_t)_{t\in[t_0,t_1)}$.
\item Let $G:[t_0,t_1)\times\R^N\times\Omega\rightarrow\R$ measurable,  bounded and such that, for all $x\in\R^N$ the process $(G(t,x,\cdot))_{t\in[t_0,t_1)}$ is adapted to the filtration $(\FR^M_t)$. Then we get, $ P$-a.s., for all $t\in[t_0,t_1)$,
\[ E[G(t,X)|\FR^M_t]=\int_{\R^N}G(t,x)dM_t(x).\]
\end{itemize}
2) Conversely, let $(\FR_t)_{t\in[t_0,t_1)}$ be a right continuous and complete filtration on $( \Omega, \FR, P)$. For all $t\in[t_0,t_1)$ and $B\in\BR(\R^n)$, set
\[ M_t(B)= P[X\in B|\FR_t].\]
Then $(M_t)$ is a martingale measure, which admits a c\`adl\`ag modification, i.e. a modification $ M$ such that,  for all measurable and bounded function $f$, $( M_tf)$ is c\`adl\`ag. It is then clear that this modification belongs to $\MR(t_0,t_1)$. In the sequel we identify all martingale measures with their c\`adl\`ag modifications.

\end{Remarks}

\noindent We write $\MR(t_0,m)$ for $\MR(t_0,T,m)$.\\

\noindent Now let
\begin{equation}
\label{w}
W(t_0,m)=\inf_{M\in\MR(t_0,m)}E\left[\Int_{t_0}^TH(s,M_s)ds\right],
\end{equation}
where we recall that
$H(t,m)=\inf_{u\in U}\sup_{v\in V}\Int_{\R^N}\ell(x,t,u,v)dm(x)$.\\
\noindent The main result of this section is the following representation formula:

\begin{Theorem}
\label{altform}
For all $(t_0,m)\in[0,T]\times\Pi_2$, it holds that
\begin{equation}
\label{wv}
 W(t_0,m)=\bV^+(t_0,m).
 \end{equation}
\end{Theorem}

\noindent{\bf Proof:}
We shall prove the two inequalities providing (\ref{wv}).\\

\noindent 1) \underline{$\bV^+(t_0,m)\leq W(t_0,m)$}\\

\noindent Fix $M\in\MR(t_0,m)$ and let $\bar u(t,m)$ be a measurable selection of
$\mbox{Argmin}_{u\in U}\sup_{v\in V}\Int_{\R^N}\ell(x,t,u,v)dm(x)$.\\
Consider the control process $\tilde u=\left(u(t,M_t)\right)_{t\in[t_0,T]}$. Remark that $\tilde u$ is adapted to the filtration $(\FR^M_t)$. Since $X$ is a realization of the probability measure $m$, the joint distribution of $(\tilde u,X)$ on $\UR(t_0)\times\R^N$ can be desintegrated as follows:
\[  P[(\tilde u,X)\in dudx]=\tilde P_x(du)m(dx),\]
for some measurable family of probability measures $(\tilde P_x,x\in\R^N)$.\\
Let now be a strategy $\alpha_M$ for Player I such that, for all $x\in\R^N$, for all measurable bounded $f:\UR(t_0)\rightarrow\R$ and for all $v\in\VR(t_0)$,
\[ E_1[f(\alpha_M(x,v))]=\int_{\UR(t_0)}f(u)\tilde P_x(du).\]
Then, for any arbitrary  $\beta\in\BR(t_0)$, strategy for Player II, we got
\[ \begin{array}{rl}
\JR(t_0,m,\alpha_M,\beta)=&
\Int_{\R^N}\int_{\UR(t_0)}\Int_{t_0}^T\ell(s,x,u(s),\beta(u)_s)ds\tilde P_x(du)dm(x) \\
\\
=&
E\left[\Int_{t_0}^T\ell(s,X,\tilde u_s,\beta(\tilde u)_s)ds\right]\\
\\
=& E\left[\Int_{t_0}^TE\left[\ell(s,X,\tilde u_s,\beta(\tilde u)_s)|\FR^M_s\right]ds\right]\\
\\
=&E\left[\Int_{t_0}^T\Int_{\R^N}\ell(s,x,\tilde u_s,\beta(\tilde u)_s)M_s(dx)ds\right].
\end{array}\]
It follows that
\[ \JR(t_0,m,\alpha_M,\beta)\leq E\left[\Int_{t_0}^T\sup_{v\in V}\Int_{\R^N}\ell(s,x,\tilde u_s,v)M_s(dx)ds\right].\]
Due to the arbitrariness of $\beta\in\BR(t_0)$ and the definition of $\tilde u$, it follows that
\[\begin{array}{rl}
 \sup_{\beta\in\BR(t_0)}\JR(t_0,m,\alpha_M,\beta)\leq &E\left[\Int_{t_0}^T\sup_{v\in V}\Int_{\R^N}\ell(s,x,\tilde u_s,v)M_s(dx)ds\right]\\
 \\
 =&E\left[\Int_{t_0}^TH(s,M_s)ds\right],
 \end{array}\]
The result follows.\\

\noindent 2) \underline{$W(t_0,m)\leq \bV^+(t_0,m)$}\\

2.0) Fix $\epsilon>0$. Let $\alpha^\epsilon\in\tilde\AR_X(t_0)$ be $\epsilon$-optimal for $\bV^+(t_0,m)$.
The assertion follows as soon as we find some $M\in\MR(t_0,m)$ and $\beta\in\BR(t_0)$ such that
\[ E\left[\Int_{t_0}^T\ell\left(X,s,\Phi(\alpha^\epsilon(X),\beta)(s)\right)ds\right]\geq E\left[\Int_{t_0}^TH(s,M_s)ds\right]-\epsilon.\]
More precisely, we shall prove by induction, that there exists $M\in\MR(t_0,m)$ and $\beta\in\BR(t_0)$, such that, if  $t_0<t_1<\ldots<t_n=T$ denotes the time grid associated to the strategy $\alpha^\epsilon$, it holds, 
for all $k\in\{ 1,\ldots,n\}$, that
\begin{equation}
\label{WV}
 E\left[\Int_{t_{k-1}}^{t_k}\ell\left(X,s,\Phi(\alpha^\epsilon(X),\beta)(s)\right)ds\right]\geq E\left[\Int_{t_{k-1}}^{t_k}H(s,M_s)ds\right]-\epsilon/n.
 \end{equation}
In parallel to the strategy  $\beta$ we shall construct a random process $Z:\Omega\times[t_0,T]\rightarrow U$ wich satisfies
\begin{equation}
\label{U}
\Phi(\alpha^\epsilon(X),\beta)(s)=(Z_s,\beta(Z)(s))\; \lambda\times P \mbox{-a.s.}.
\end{equation}

\noindent 2.1) Because of the delay, on $[t_0,t_1)$, the control of Player I doesn't depend on the strategy of Player II. Since, by definition, $\alpha^\epsilon$ is a measurable map from $\Omega_1\times\R^N$ to $\AR(t_0)$, this control depends also in a measurable way of $(\omega_1,x)$. More precisely, we can set, for any $(x,\omega,v)\in\R^N\times\Omega\times \VR(t_0)$, $u^\epsilon(x,\omega)=\alpha(x,\omega_1,v)|_{[t_0,t_1)}\in\UR(t_0,t_1)$ with $\omega=(\omega_0,\omega_1)$. 
We finally define the random process $(Z_s)$ on the time interval $[t_0,t_1)$ by:  $Z_s=u^\epsilon(X)_s, s\in[t_0,t_1)$.\\
Then, for all $\beta\in\BR(t_0)$ and $x\in\R^N$, a trivial application of the fixed point relation (\ref{fix}) of Lemma \ref{fixpoint} shows that, on the same interval $[t_0,t_1)$, (\ref{U}) is satisfied.\\
We denote by $(\FR_s,s\in[t_0,t_1))$  the completed right-continuous filtration generated by the process $(Z_s)$. 
Let $M\in\MR(t_0,m)$ be the martingale measure defined on $[t_0,t_1)$ by: for all $B\in\BR(\R^N)$, 
\[
 M_s(B)=P[X\in B|\FR_s],\;  s\in [t_0,t_1).
\]

\noindent 2.2) Suppose now that,
 for some $k\in\{ 0,\ldots,n-1\}$,  for all $u\in\UR(t_0)$, $\beta(u)$ is defined on $[t_0,t_k)$, as well as $(Z_s)$ satisfying (\ref{U}) on $[t_0,t_{k+1})$. Let $M\in\MR(t_0,m)$ still be a c\`adl\`ag version of the martingale measure $M^0$ defined on $[t_0,t_k+1)$ by 
\[
 \Int_{\R^N}f(x)M^0_s(dx)=E[f(X)|\FR_s],\; s\in [t_0,t_{k+1}).
\]
 We suppose also that the inequality (\ref{WV}) holds for all $j\in\{1,\ldots,k\}$.\\
Now we extend the strategy $\beta$  on $[t_k,t_{k+1})$:\\
First of all, for $u\in U$ and $s\in[t_k,t_{k+1})$, we denote by $\bar v(s,u,m)$ a measurable selection of $\mbox{Argmax}_{v\in V}\Int_{\R^N}\ell(x,s,u,v)m(dx)$.\\
Remark that, for all $s\in[t_0,t_{k+1})$, $M_{s-}$ is measurable with respect to $\FR_{s-}=\sigma(Z_r,r< s)\vee\NR$, where $\NR$ is the set of null-sets for the probability $P$.\\
Therefore there exists some measurable application $\Psi:[t_0,t_{k+1})\times\UR(t_0)$, such that, for all  $u\in\UR(t_0)$  and all fixed $s\in[t_0,t_{k+1})$, $ \Psi(s,u)$ depends only  on the restriction of $u$ on $[t_0,s)$ : $\Psi(s,u)=\Psi(s,(u_r,r\in[t_0,s))$ and such that
$\bar v(s,Z_s,M_{s-})=\Psi(s,(Z_r,r\in[t_0,s]))$ $\lambda\times P$-a.s..  Now we set, for any $(u(s), s\in[t_0,T])\in\UR(t_0)$ and $s\in[t_k,t_{k+1})$,
$\beta^0(u)(s)=\Psi(s,(u_r,r\in[t_0,s]))$.\\ 
By Lemma \ref{delaybeta}, there exists, for all $u\in\UR(t_0)$ some delay $\delta(u)>0$ such that
\[ \left|E\Int_{t_k}^{t_{k+1}}\ell(X,s,u_s,\beta^0(u)_s)ds
- E\Int_{t_k}^{t_{k+1}}\ell(X,s,u_s,\beta^0(u)_{(s-\delta(u))_+})ds\right|\leq \epsilon/n.\] 
We can extend now the strategy $\beta$ on $[t_k,t_{k+1})$ by setting
\[ \forall u\in\UR(t_0), \forall s\in[t_k,t_{k+1}),\; \beta(u)_s=\beta^0(u)_{(s-\delta(u))_+}.\]
We get then
\[
\begin{array}{rl}
E\left[\Int_{t_k}^{t_{k+1}}\ell(X,s,\Phi(\alpha^\epsilon(X),\beta)(s))ds\right]
=& E\left[\Int_{t_k}^{t_{k+1}}\ell(X,s,Z_s,\beta^0(Z)(s-\delta(Z))_+)ds\right]\\
\\
\geq &E\left[\Int_{t_k}^{t_{k+1}}\ell(X,s,Z_s,\beta^0(Z)(s))ds\right]-\epsilon/n.
\end{array}
\]
Further, since the filtration $(\FR_s)$ is generated by $(Z_s)$ and by the definition of $M$, it holds, for all $s\in[t_{k-1},t_k)$,
\[ E\left[\ell(X,s,Z_s,\beta^0(Z)(s))|\FR_s\right]
=\int_{\R^N}\ell(x,s,Z_s,\beta^0(Z)(s))dM_s(x)\; P\mbox{-a.s.}.\]
Integrating over $[t_k,t_{k+1})$ and taking the expectation of both sides, we get
\[\begin{array}{rl} 
E\left[\Int_{t_k}^{t_{k+1}}\ell(X,s,Z_s,\beta^0(Z)(s))ds\right ]=&
 E\left[\Int_{t_k}^{t_{k+1}}\int_{\R^N}\ell(x,s,Z_s,\bar v(s,Z_s,M_{s-}))dM_s(x)ds\right ]\\
\\
=&
 E\left[\Int_{t_k}^{t_{k+1}}\int_{\R^N}\ell(x,s,Z_s,\bar v(s,Z_s,M_s))dM_s(x)ds\right ]\\
\\
=&E\left[\Int_{t_k}^{t_{k+1}}\left(\sup_{v\in V}\Int_{\R^N}\ell(x,s,Z_s,v)M_s(dx)\right)ds\right]
\\
\geq & E\left[\Int_{t_k}^{t_{k+1}}H(s,M_s)ds\right].
\end{array}\]

It remains now to extend the process $Z$ on $[t_{k+1},t_{k+2})$ : Using  the fix point Lemma \ref{fixpoint}, we define first the map $u^\epsilon:\R^N\times\Omega\rightarrow\UR(t_{k+1},t_{k+2})$   by
\[ \forall (x,\omega_1)\in\R^N\times\Omega_1, \Phi(\alpha^\epsilon(x,\omega_1),\beta)=(u^\epsilon(x,\omega_1),\beta(u^\epsilon(x,\omega_1))\]
 and set then, for all $s\in[t_{k+1},t_{k+2})$,
\[ Z_s(\omega)=u^\epsilon(X(\omega_0),\omega_1).\]
The result follows.
\QED

In the remaining part of this section, we prove that the infimum in (\ref{w}) is attained up to  a change of probability space. Let  ${\mathcal P}$ be the set of Radon nonnegative measures $m$ on $[-1,T+1]\times \R^N$, with total mass $T+2$, with a first marginal equal to the Lebesgue measure on $[-1,T+1]$ and such that 
$$
\int_{(-1, T+1)\times \R^N} |x|^2 m(dt,dx) <+\infty\;.
$$
We endow  ${\mathcal P}$ with the weak topology. Any measure $m\in {\mathcal P}$ can be desintegrated into $m(dx,ds)= m_s(dx) \ ds$, where $m_s$ is a probability measure on $\R^N$ for almost every $s$. In the sequel we systematically identify measure $m\in {\mathcal P}$ and Borel families of probability measures $(m_s)_{s\in [-1,T+1]}$. 
We denote by $\Gamma$ the subset of ${\mathcal P}$ consisting in measures $m\in {\mathcal P}$ which have a c\`adl\`ag representative for the ${\bf d}_1$ distance. The topology on  $\Gamma$ is the one induced by ${\mathcal P}$. We denote by $t\to m_t$ the canonical process on $\Gamma$. For a fixed initial condition $(\bar t, \bar m)\in [0,T]\times \Pi_2$ we finally  denote by $\overline {\mathcal M}(\bar t,\bar m)$ the set of Borel probability measures $M$ on $\Gamma$ under which, 
\begin{itemize}
\item[(i)] $m_{\bar t}= \bar m$ for $t<\bar t$, 
\item[(ii)] $m_t \in {\mathcal D}$  for $t\geq T$, where ${\mathcal D}$ is the set of single Dirac masses on $\R^N$, 
\item[(iii)] for any continuous, bounded map $\varphi:\R^N\to \R$,  the  process $t\to \int_{\R^N} \varphi(x)\ dm_t(x)$ is a c\`adl\`ag martingale for the filtration $({\mathcal F}_t)$ generated by the process $(m_t)$. 
\end{itemize}
Note that the law of any martingale measure $M\in \MR(t_0,\bar m)$ belongs to $\overline {\mathcal M}(\bar t,\bar m)$.

\begin{Theorem}\label{theo:existsoptimal}
For any $(\bar t,\bar m)\in[0,T]\times \Pi_2$, there is some $\bar M\in \overline {\mathcal M}(\bar t,\bar m)$ such that
$$
\bV^+(\bar t,\bar m)= \int_{\Gamma} \Int_{\bar t}^TH(s,m_s)\ ds\ d\bar M(m)\;.
$$
\end{Theorem}

Theorem \ref{theo:existsoptimal} is an easy consequence of  the following compactness property of the set $\overline {\mathcal M}(\bar t,\bar m)$: 
\begin{Proposition}\label{prop:compactM} Let $(\bar t,\bar  m)\in [0,T]\times \Pi_2$. 
Let $(M^n)$ be a sequence in $\overline {\mathcal M}(\bar t,\bar  m)$. Then there is a subsequence $(M^{n'})$ and a set $I$ of full measure in $[-1,T+1]$ such that $(M^{n'})$  converges weakly (as a measure over ${\mathcal P}$) to some $M\in \overline {\mathcal M}(\bar t,\bar  m)$ and  $(M^{n'}_t)$ converges weakly (as a measure over $\Delta(\R^N)$) to $M_t$ for any $t\in I$. 
\end{Proposition}

The proof of Proposition \ref{prop:compactM} is a variation on a paper by Meyer and Zheng \cite{MeyerZheng}. Let us start with some sufficient condition for a measure on ${\mathcal P}$ to be supported by $\Gamma$. For this, we need some notations. For $m\in {\mathcal P}$ and $\varphi\in {\mathcal C}^0$ with at most a linear growth, we denote by $t\to m_t(\varphi)$ the real-valued measurable map $t\to \int_{\R^N} \varphi(x)dm_t(x)$. 


 Let $t\to f_t$ on $[-1,T+1]$ be a real valued measurable map and $u<v$. Following \cite{MeyerZheng}, 
we define  the number of crossing $N^{uv}(f)$ as the largest number $k$ for which we can find times $-1\leq  t_1<t_1'<\dots <t_k<t_k' \leq T+1$ such that each set of the form  $(-1, t_1)\cap  f^{-1}((-\infty, u))$, $(t_i,t_{i}')\cap f^{-1}( (v,+\infty))$, $(t_i',t_{i+1})\cap f^{-1}((-\infty,u))$ and  $(t_k',T+1)\cap f^{-1}((v,+\infty))$ has a positive measure. Let $Lip_1$ be the set of $1-$Lipschitz continuous maps $\varphi$ on $\R^N$ such that $\varphi(0)=0$.

\begin{Lemma}\label{lem:charGamma} Let $M$ be a probability measure on ${\mathcal P}$ such that, for some constant $C$,   
\be\label{hyp1lemacharGamma}
M \left[ \sup_{t\in [-1, T+1]} \int_{\R^N} (|x|^2+1) \ dm_t(x)\; >\; R\right ]\leq \frac{C}{R}\qquad \forall R\geq 1
\ee
and
$$
\int_{{\mathcal P}}  N^{uv}(m(\varphi))\ dM(m) \leq \frac{|u|+C}{v-u}\qquad \forall u<v, \; \varphi\in Lip_1\;.
$$
Then  $M(\Gamma)=1$.
\end{Lemma}

\noindent {\bf Proof of the Lemma: }  We endow $Lip_1$ with the distance
$$
d(\varphi, \psi) =  \left\|\frac{\varphi-\psi}{ |x|^2+1} \right\|_{L^\infty(\R^n)}\;.
$$
Then $Lip_1$ is a compact set. Fix   $\ep>0$, $R\in (0, 1/(2\ep))$, $u<v$ such that $v-u>2R\ep$ and $\varphi^1, \dots, \varphi^n$ in $Lip_1$ such that $\bigcup_{k=1}^n B(\varphi^k,\ep) \supset Lip_1$, where $B(\varphi^k,\ep)$ is the closed ball centered at $\varphi^k$ and of radius $\ep$ in $Lip_1$.  In view of our assumption, 
$$
M\left[ N^{(u+R\ep)(v-R\ep)}(m(\varphi^k))\geq \frac{n}{\ep} \right]\leq \frac{|u|+R\ep+ C}{v-u-2R\ep}\ \frac{\ep}{n}
$$
So 
$$
M\left[ \sup_{k\in \{1,\dots, n\}} N^{(u+\ep)(v-\ep)}(m(\varphi^k))\geq \frac{n}{\ep} \right]\leq \frac{|u|+R\ep+ C}{v-u-2R\ep}\ \ep
$$
In view of assumption (\ref{hyp1lemacharGamma}), we have
$$
M \left[ \sup_{t\in [-1,T+1]} \int_{\R^N} (|x|^2+1) \ dm_t(x) \; > \; R\right] < \frac{C}{R}\;.
$$
Let $E_R$ denote the set of $m\in {\mathcal P}$ such that $\ds \sup_{t\in [-1,T+1]} \int_{\R^N} (|x|^2+1) \ dm_t(x) \; \leq \; R$. 
Let $\varphi\in Lip_1$ and $m\in E_R$. By definition of the $(\varphi^k)$,  there is some $k\in \{1,\dots, n\}$ such that $d(\varphi,\varphi^k)\leq \ep$. Then 
$$
\sup_{t\in [-1, T+1]} |m_t(\varphi)-m_t(\varphi^k)| \leq \sup_{t\in [-1, T+1]} (|x|^2+1)dm_t(x) \ d(\varphi,\varphi^k) \leq 
R\ep\;.
$$
Hence 
$$
\sup_{\varphi\in Lip_1} N^{uv}(m(\varphi)) \leq \sup_{k\in \{1,\dots, n\}} N^{(u+R\ep)(v-R\ep)}(m(\varphi^k))\;, 
$$
so that 
$$
M\left[ \left\{\sup_{\varphi\in Lip_1} N^{uv}(m(\varphi)) \geq \frac{n}{\ep}\right\}\; \cup E_R^c \right]\leq \frac{|u|+R\ep+ 1}{v-u-2R\ep}\ \ep+\frac{C}{R}\;.
$$
Letting $\ep\to 0$ and then $R\to +\infty$ we obtain, for any $u<v$ and $M-{\rm a.s.}$
\be\label{Nuv}
\sup_{\varphi\in Lip_1} N^{uv}(m(\varphi)) <+\infty \qquad{\rm and }\qquad 
 \sup_{t\in [-1, T+1]} \int_{\R^N} (|x|^2+1) \ dm_t(x)<+\infty \;.
\ee
It remains to show that any $m$ satisfying the above inequalities coincides (up to a subset of measure $0$) with an element of $\Gamma$. 
 
Let us fix $t\in (0,T]$ and let  $L_t(m)$ the left essential upper limit of $(m_s)$ as $s\to t^-$: 
$$
\nu\in L_t(m) \; \Leftrightarrow \; \forall \delta>0, \; \left| \{ s\in (t-\delta,t), \; m_s\in B(\nu,\delta)\}\right| >0  
$$
where $B(\nu,\delta)$ is the ball centered at $\nu$ and of radius $\delta$ for the Monge-Kantorovich $ {\bf d}_1$ distance. 
We claim that $L_t(m)$ is a singleton. Indeed, note first that $L_t(m)$ is not empty by compactness. Assume that $L_t(m)$ is not a singleton. Then there would exist some $\nu_1\neq \nu_2$ in $L_t(m)$. Since
$$
\sup_{\varphi\in Lip_1} \int_{\R^N} \varphi d(\nu_1-\nu_2)= {\bf d}_1(\nu_1,\nu_2)\;, 
$$
there is some $\varphi\in Lip_1$ such that $ \int_{\R^N} \varphi d(\nu_1-\nu_2)=\eta>0$. Let us set $u=  \int_{\R^N} \varphi d\nu_2+ \eta/4$ and 
$v= \int_{\R^N} \varphi d\nu_1-\eta/4$. Then $v>u$ and, by definition of $L_t(m)$, $N^{uv}(m(\varphi))=+\infty$. This contradicts (\ref{Nuv}). 
Therefore $L_t(m)$ is a singleton for any $t$. In the same way, one can prove that  the right essential upper limit $R_t(m)$ of $(m_s)$ as $s\to t^+$
is reduced to a singleton for any $t\in [0,T)$. One then easily checks that $t\to R_t(m)$ is right-continuous, has a left-limit on $(0,T)$ and that $m_t=R_t(m)$ for a.e. $t$. Therefore $m$ belongs to $\Gamma$ because $m$ has a c\`adl\`ag representative $R_\cdot(m)$. 
%
\QED

\noindent {\bf Proof of Proposition \ref{prop:compactM}: } Let $(M^n)$ be a sequence in $\overline {\mathcal M}$. 
Since, for any $k>0$ and for any fixed $n$, the random process $t\to \int_{\R^N} (|x|^2\wedge k)dm_t(x)$ is a martingale under $M^n$, we have, from Doob's maximal inequality,
$$
M^n\left[\sup_{t\in [-1, T+1]} \int_{\R^N} (|x|^2\wedge k)dm_t(x)\geq R \right] \leq \frac{1}{R} \int_{\R^N} |x|^2 d\bar m(x)=: \frac{\bar C}{R}
\qquad \forall R\geq 1\;.
$$
Letting $k\to +\infty$, we obtain 
$$
M^n\left[\sup_{t\in [-1, T+1]} \int_{\R^N} |x|^2\ dm_t(x)\geq R \right] \leq \frac{\bar C}{R}
\qquad \forall R\geq 1\;.
$$
Since the set $\ds \left\{m\in {\mathcal P},\; \sup_{t\in [-1, T+1]} \int_{\R^N} |x|^2\ dm_t(x)\leq R\right\}$ is compact, the sequence $M^n$  is tight and therefore  a subsequence of $(M^n)$, still denoted $(M^n)$, converges weakly  on ${\mathcal P}$ to some probability law $M$. 
Note that 
$$
M\left[\sup_{t\in [-1, T+1]} \int_{\R^N} |x|^2\ dm_t(x)> R \right] \leq \frac{\bar C}{R}
\qquad \forall R\geq 1\;.
$$
Let us now show that $M$ is supported by $\Gamma$. Fix this, let us denote by  $\E^n$ and $\E$ the expectations with respect to $M^n$ and $M$ respectively. 
Fix $\varphi\in Lip_1\cap {\mathcal C}^0_b(\R^N)$, where ${\mathcal C}^0_b(\R^N)$ denotes the set of continuous, bounded function in $\R^N$. Our aim is to estimate $\E\left[N^{uv}(m (\varphi))\right] $. Since, under $M^n$,  the  process $t\to \int_{\R^N} \varphi(x)\ dm_t(x)$ is a martingale, its conditional variation (in the sense of \cite{MeyerZheng}) can be estimated by
$$
\sup_{t\in [-1, T+1]} \E^n\left[ \left| \int_{\R^N} \varphi(x)\ dm_t(x)\right|\right]
\leq   \sup_{t\in [-1, T+1]}  \E^n\left[ \int_{\R^N} (|x|^2+1)\ dm_t(x)\right] \ \left\|\frac{\varphi}{|x|^2+1}\right\|_\infty
\leq \bar C+1
$$
So, by 
Lemma 3 of \cite{MeyerZheng}, we have, for any $u<v$, 
$$
\E^n\left[ N^{uv}(m(\varphi))\right] \leq \frac{|u|+ \bar C+1}{v-u}\;.
$$
Since $N^{uv}$ is lower semicontinuous for the weak topology, we get
$$
\E\left[N^{uv}(m (\varphi))\right] \leq \frac{|u|+ \bar C+1}{v-u} \;.
$$
Now we can relax the boundedness assumption on $\varphi$ to get that the above inequality holds for any $\varphi\in Lip_1$. 
Then Lemma \ref{lem:charGamma} implies that $M(\Gamma)=1$. Note also that, under $M$, $m_{0^-}=\bar m$ and $m_t\in {\mathcal D}$ for any $t\geq T$ because ${\mathcal D}$ is closed under the weak topology. \\

We now prove the finite dimensional convergence. Since ${\mathcal P}$ is a Polish space and $M^n$ converges weakly to $M$, Skorokhod's theorem implies that there is a probability space $(\Omega, {\mathcal A}, \P)$ and random variables ${\bf M}^n$ and ${\bf M}$ on ${\mathcal P}$ with respective law $M^n$ and $M$ and such that ${\bf M}^n$ almost surely converges to ${\bf M}$. Note that the processes $t\to {\bf M}^n_t$ and $t\to {\bf M}_t$ are c\`adl\`ag $\P-$a.s. because $M^n$ is supported by $\Gamma$. Let now $(\varphi^k)$ be a sequence which is dense in ${\mathcal C}^0_b(\R^N)$ (for the uniform convergence on compact subsets of $\R^N$). From our assumption, the process  $t\to {\bf M}^n_t\varphi^k$ is a martingale. So, by a diagonal argument, (the proof of) Theorem 5 of \cite{MeyerZheng} gives the existence of a set of full measure $I$ in $[0,T]$ and a subsequence (again denoted $({\bf M}^n)$) such that $({\bf M}^n_t\varphi^k)$ converges to ${\bf M}_t\varphi^k$ as $n\to +\infty$ $\P-$a.s. Since $(\varphi^k)$ is dense, this easily implies that $({\bf M}^n_t)$ weakly converges to ${\bf M}_t$ for any $t\in I$  $\P-$a.s.. This means that $M^n_t$ converges weakly to $M_t$ (as measures on $\Delta(\R^N)$) for any $t\in I$: for any $\varphi\in {\mathcal C}^0_b(\Delta(\R^N))$, 
$$
\int_{\Delta(\R^N)} \varphi(\nu)dM^n_t(\nu)= \E\left[ \varphi({\bf M}^n_t)\right]\to \E\left[ \varphi({\bf M}_t)\right]=\int_{\Delta(\R^N)} \varphi(\nu)dM_t(\nu)\qquad \forall t\in I\;.
$$

Let us finally show that $M$ belongs to $\overline {\mathcal M}$. For this it only remains to show that, for any continuous bounded map $\varphi:\R^N\to \R$, the process ${\bf M}_t\varphi$ is a martingale for the filtration generated by ${\bf M}_t$. Let  $k_1, \dots, k_m$ be  indices, $t_1\leq \dots\leq  t_m, t$ be times in $I$ and $f_1,\dots, f_m, f$ be bounded continuous maps. Then, for any $t\in [0,T]$, and since ${\bf M}^n_t\varphi$ is a martingale,  
$$
\E\left[ f_1({\bf M}^n_{t_1}\varphi^{k_1})\dots f_m({\bf M}^n_{t_m}\varphi^{k_1})	 f({\bf M}^n_t\varphi)\right]
= 
\E\left[ f_1({\bf M}^n_{t_1}\varphi^{k_1})\dots f_m({\bf M}^n_{t_m}\varphi^{k_1})	 f({\bf M}^n_{t\wedge t_m}\varphi)\right]
$$
Letting $n\to +\infty$ and using the finite dimensional convergence, we get
$$
\E\left[ f_1({\bf M}_{t_1}\varphi^{k_1})\dots f_m({\bf M}_{t_m}\varphi^{k_1})	 f({\bf M}_t\varphi)\right]
= 
\E\left[ f_1({\bf M}_{t_1}\varphi^{k_1})\dots f_m({\bf M}_{t_m}\varphi^{k_1})	 f({\bf M}_{t\wedge t_m}\varphi)\right]
$$
Since ${\bf M}$ is c\`adl\`ag, the $\sigma-$algebra generated by the family of random variables $f_1({\bf M}_{t_1}\varphi^{k_1})$, $\dots$, $f_m({\bf M}_{t_m}\varphi^{k_1})$ (where $k_1, \dots, k_m$, $t_1\leq \dots\leq  t_m$ and $f_1,\dots, f_m$ are as above) is equal to $(\sigma({\bf M}_s, \; s\leq t))_{t\in [0,T]}$. Therefore ${\bf M}\varphi$ is a martingale and the proof of the Proposition is complete.
\QED

\noindent {\bf Proof of Theorem \ref{theo:existsoptimal}: } In view of Theorem \ref{altform} we know that there exists a sequence ${\bf M}^n\in\MR(\bar t,\bar m)$ such that 
$$
\lim_{n\to +\infty} E\left[\Int_{\bar t}^TH(s,{\bf M}^n_s)ds\right]= \V^+(\bar t, \bar m)\;.
$$
Let $M^n$ be the law of ${\bf M}^n$ and recall that $M^n\in \overline {\mathcal M}(\bar t,\bar  m)$ and that 
$$
E\left[\Int_{\bar t}^TH(s,{\bf M}^n_s)ds\right]= \int_{\Gamma}\int_{\bar t}^T H(s,m_s)\ ds\ dM^n(m)\;.
$$
According to Proposition \ref{prop:compactM} there is a measure $\bar M\in  \overline {\mathcal M}(\bar t,\bar  m)$, a subsequence of the $(M^n)$, still denoted $(M^n)$, and a set $I$ of full measure in $[0,+\infty)$ such that  $M^{n}_t$ converges weakly (as a measure over $\Delta(\R^N)$) to $\bar M_t$ for any $t\in I$.  Then, for any $t\in I$,  
$$
\lim_{n\to+\infty} 
\Int_\Gamma H(t,m_t)\ d M^n(m)=
\Int_\Gamma H(t,m_t)\ d \bar M(m)
$$
because the map $\mu\to H(t,\mu)$ is continuous in $\Pi_2$ for the ${\bf d}_1$ distance. Then,  since $H$ is bounded, we conclude by Lebesgue convergence Theorem that 
$$
  \V^+(\bar t, \bar m)= \int_{\bar t}^T \Int_\Gamma H(t,m_t)\ d \bar M(m)\;.
$$

\QED

\noindent A simple application of the above characterization of the value function is the following  dynamic programming principle: 
\begin{Corollary}\label{PGD} We have
\begin{equation}
\label{dpp}
\bV^+(t,m)= \inf_{M\in\MR(t,m)}E\left[\Int_{t}^{t+h}H(s,M_s)ds+ \bV^+(t+h, M_{t+h})\right]\qquad \forall h\in [0,T-t]\;.
\end{equation}

\noindent Moreover there exists an enlargement of $(\Omega,\FR,P)$, and, on it, some $\tilde M$ as above such that 
$$
\bV^+(t,m)= E\left[\Int_{t}^{t+h}H(s,\tilde M_s)ds+ \bV^+(t+h, \tilde M_{t+h})\right]\qquad \forall h\in [0,T-t]\;,
$$
\end{Corollary}

\noindent{\bf Proof:} As usual, we split the relation (\ref{dpp}) in two inequalities which we prove separately.
 Let us denote by $\bV^h(t,m)$ the right hand side of (\ref{dpp}) and, for $\epsilon>0$, let $M^\epsilon\in\MR_f(t,m)$ be an $\epsilon$-optimal martingale measure for $\bV^h(t,m)$ with finite support (i.e., the support of $M^\epsilon$ is concentrated on a finite number of measures).\\
Let $A_1,\ldots,A_n$ be the atoms af the $\sigma$-algebra $\FR^{M^\epsilon}_{t+h}$ and $m_1,\ldots, m_n\in \Pi_2$ the values taken by $M^\epsilon_{t+h}$ on these sets. Remark that, for all $i\in\{ 1,\ldots,n\}$, $m_i$ is the law of the restriction of $X$ on $A_i$.\\
For each $i\in\{ 1,\ldots, n\}$, let $M^{(i)}\in\MR(t+h,m_i)$ be $\epsilon$-optimal for $V(t+h,m_i)$ and define some measure valued processes $(M^i_s,s\in [t+h,T])$ such that $M^i|_{A_i\times[t+h,T]}\stackrel{(d)}=M^{(i)}$
and with terminal value $M_T=\delta_X$.\\
Finally we set 
\[ \bar M^\epsilon_s=\left\{ \begin{array}{ll}
M^\epsilon_s&\mbox{ for } s<t+h,\\
M^i_s&\mbox{ on } A_i\times[t+h,T], \; i\in\{ 1,\ldots,n\}.
\end{array}\right.\]
It is easy to prove that $\bar M^\epsilon\in\MR(t,m)$. Therefore it holds that
\[\begin{array}{rl}
\bV^+(t+h,\bar M_{t+h})=&
\sum_{i=1}^n\ind_{A_i}\bV^+(t+h,m_i)\\
\geq & \sum_{i=1}^n\ind_{A_i} E\left[\int_{t+h}^TH(s,M^{(i)}_s)ds\right]-\epsilon\\
=& \sum_{i=1}^n\ind_{A_i}E\left[\int_{t+h}^TH(s,M^{i}_s)ds|A_i\right]-\epsilon\\
= &
E\left[\int_{t+h}^TH(s,\bar M^\epsilon_s)ds|\FR^{M^\epsilon}_{t+h}\right]-\epsilon.
\end{array}\]
It follows that 
\[ \begin{array}{rl}
\bV^h(t,m) \geq &E\left[ \int_t^{t+h}H(s, M^\epsilon_s)ds\right] +E\left[\int_{t+h}^TH(s,\bar M^{\epsilon}_s)ds\right]-2\epsilon\\
=&E\left[\int_t^TH(s,\bar M^\epsilon_s)ds\right]-2\epsilon\\
\geq & \bV^+(t,m)-2\epsilon.
\end{array}\]
Since the last relation holds true for all $\epsilon>0$, we get our first inequality:
\[ \bV^h(t,m) \geq \bV^+(t,m).\]

\noindent 1.2 Consider now $M^\epsilon\in\MR_f(t,m)$ being $\epsilon$-optimal for $\bV^+(t,m)$. 
Again we denote by $A_1,\ldots, A_n$ the atoms of $\FR^{M^\epsilon}_{t+h}$ and by $m_1,\ldots, m_n$ the values taken by $M^\epsilon_{t+h}$.\\
Then we can write
\[ \bV^+(t,m)\geq E\left[\int_t^{t+h}H(s,M^\epsilon_s)ds+\sum_{i=1}^n\ind_{A_i}E\left[\int_{t+h}^TH(s,M^\epsilon_s)ds|A_i\right]\right]-\epsilon.\]
Now, restricted on $A_i$, $M^\epsilon$ is a martingale measure with initial condition $m_i$. Therefore it holds that
\[ E\left[\int_{t+h}^TH(s,M^\epsilon_s)ds|A_i\right]\geq \bV^+(t+h,m_i).\]
The dynamic programming principle follows. \\

\noindent 2. By Theorem \ref{theo:existsoptimal}, we can enlarge the probability space $(\Omega,{\mathcal F}, P)$ such that there exists $\tilde M\in\MR(t,m)$ which is optimal for $\bV^+(t,m)$.
For $\epsilon>0$, let $M^\epsilon\in\MR_f(t,m)$ such that $E[\int_{t+h}^T{\bf d}_1(M^\epsilon_s,\tilde M_s)ds]\leq \epsilon$.  Since the trajectories of  $\tilde M$ are c\`adl\`ag, we can choose $M^\epsilon$ such that 
$E[{\bf d}_1(M^\epsilon_{t+h},\tilde M_{t+h})]\leq \epsilon$.
Using the fact that $V$ is Lipschitz in $m$, we get now
\[\begin{array}{rl}
 E[\bV^+(t+h,\tilde M_{t+h})]\leq & E[ \bV^+(t+h,M^\epsilon_{t+h})]+ C\epsilon\\
=&E[\sum_i\bV^+(t+h,m_i)\ind_{A_i}]+C\epsilon,
\end{array}\]
where $A_1,\ldots, A_n$ are the atoms of $\FR^{M^\epsilon}_{t+h}$ and $m_1,\ldots, m_n$ the values taken by $M^\epsilon_{t+h}$ on them.\\
Now, since, for each $i\in\{ 1,\ldots,n\}$, the restriction of $(M^\epsilon_s)_{t+h\leq s\leq T}$ on $A_i$ belongs to $\MR(t+h,m_i)$, and then, using the fact that $H$ is Lipschitz, we have
\[\begin{array}{rl}
E[\sum_i\bV^+(t+h,m_i)\ind_{A_i}]\leq  & E[\int_{t+h}^TH(s,M^\epsilon_s)ds]\\
\leq & E[\int_{t+h}^TH(s,\tilde M_s)ds]+C\epsilon.\\
\end{array}\]
It follows that
\begin{equation}
\label{Vm}
  E[\bV^+(t+h,\tilde M_{t+h})]\leq E[\int_{t+h}^TH(s,\tilde M_s)ds].
\end{equation}
Finally, using the first part of the Lemma, we get
\[\begin{array}{rl}
E\left[\int_t^{t+h}H(s,\tilde M_s)ds+\bV^+(t+h,\tilde M_{t+h})\right]\leq &
E\left[\int_t^{t+h}H(s,\tilde M_s)ds+\int_{t+h}^TH(s,\tilde M_s)ds\right]\\
=&V(t,m)\\
=&\inf_{M\in\MR(t,m)}E\left[\int_t^{t+h}H(s,M_s)ds+\bV^+(t+h,M_{t+h})\right]\\
\leq & E\left[\int_t^{t+h}H(s,\tilde M_s)ds+\bV^+(t+h,\tilde M_{t+h})\right].
\end{array}\]
The result follows.
\QED


\section{Characterization of the value function.}

In this section, we give two different characterizations of the value function $\bV^+$ as 
a viscosity solution of some deterministic functional equations. The first one, stating that $\bV^+$ is the largest subsolution of some equation, is interesting by its simplicity. The second one is purely local, but is more involved in its formulation and derivation.

\subsection{Characterization of the upper value function as largest subsolution of some HJI-equation.}

\begin{Proposition} 
\label{largest}
The value function $\bV^+$ is the largest continuous map on $[0,T]\times\Pi_2$ which satisfies\\
i) For all $m\in\Pi_2$, $\bV^+(T,m)=0$,\\
ii) $\bV^+$ is convex with respect to $m$,\\
iii) for all $m\in\Pi_2$, $\bV^+(\cdot,m)$ is a subsolution (in the viscosity sense) of the following ordinary differential equation
\begin{equation}
\label{ode}
 \frac{d}{dt} \bV^+(t,m)+H(t,m)=0 \qquad {\rm in }\; (0,T)\ .
\end{equation}
\end{Proposition}

\noindent{\bf Proof:} The map $\bV^+$ clearly satisfies $\bV^+(T,m)=0$ for all $m\in\Pi_2$. Moreover  $\bV^+$ is convex in $m$ by Lemma \ref{convex}.
Let us now check that, for any $m\in\Pi_2$,  $\bV^+(\cdot,m)$ is a subsolution of (\ref{ode}).  For $t\in [0, T)$, let $\phi \in C^1([0,T])$ be such that $\bV^+(\cdot,m)\leq\phi$ and $\bV^+(t,m)=\phi(t)$. Let  $M\in\MR(t,m)$ be the  martingale measure defined by
\[ M_r=\left\{\begin{array}{rl}
m&\mbox{ if } r\in[t,T),\\
\delta_x & \mbox{with probability $m$ for $t\geq T$}.
\end{array}\right.\]
By  dynamic programming (Corollary \ref{PGD}), we have, for any $h\in [0,T-t]$,  
\[\begin{array}{rl}
\phi(t)=\bV^+(t,m) \leq   & E\left[ \Int_t^{t+h}H(r, M_r)dr+ \bV^+(t+h, M_{t+h})\right]
\\
\leq & h\ H(r,m)+ \phi(t+h)\ .
\end{array}\]
The result follows by dividing by $h$ and letting $h\to 0$.\\

We now prove that $\bV^+$ is the largest subsolution. 
Let $w:[0,T]\times\Pi_2\rightarrow\R$ be convex with respect to $m$, with $w(T,\cdot)=0$ and such that $w(\cdot,m)$ is a subsolution of (\ref{ode}). We have to prove that $w\leq \bV^+$. 
For $(t_0,m_0)\in[0,T]\times\Pi_2$, let $t_0<t_1<\ldots< t_n=T$ be a subdivision of $[t,T]$, with $t_{k+1}-t_k=\tau>0$.\\
For $m\in\Pi_2$, $k\in\{ 1,\ldots,n-2\}$ and $\epsilon>0$, let 
$w^\epsilon(t)=\sup\{ w(s,m)-\frac 1{2\epsilon}|t-s|^2, s\in[0,T]\}$ be the sup-convolution of $w$ on the interval $[t_k,t_{k+1}]$.
Then  (see e.g. \cite{bcd}) the map $w^\epsilon$ is a locally Lipschitz subsolution on $[t_k,t_{k+1}]$ of 
\[ (w^\epsilon)'(t)+H(t,m)\geq -o(1).\]
It follows that
\[ w^\epsilon(t_{k+1})-w^\epsilon(t_k)+\int_{t_k}^{t_{k+1}}H(s,m)ds\geq -\tau o(1).\]
Since $w^\epsilon\searrow w(\cdot,m)$ as $\epsilon\searrow 0$, we get for $w(\cdot,m)$:
\begin{equation}
\label{wm}
 w(t_{k+1},m)-w(t_k,m)+\int_{t_k}^{t_{k+1}}H(s,m)ds\geq 0.
 \end{equation}
Let  now $M\in\MR(t_0,m_0)$ and let us define $M^\tau\in\MR(t_0,m_0)$ by
$$
M^\tau_s = M_{t_k} \qquad {\rm for }\; s\in [t_k, t_{k+1})
$$ 
From (\ref{wm}), we get
\[
E\left[ w(t_{k+1},M_{t_{k+1}})-w(t_k,M_{t_{k+1}})+\int_{t_k}^{t_{k+1}}H(s,M_{t_{k+1}})ds\right]\geq 0.
\]
where, since $w$ is convex and $M$ a martingale measure, 
\[ E\left[ w(t_{k},M_{t_{k+1}})|\FR_{t_k}\right]\geq w(t_{k},M_{t_{k}}).\]
It follows that
\begin{equation}
\label{wH}
E\left[ w(t_{k+1},M_{t_{k+1}})-w(t_k,M_{t_{k}})+ \Int_{t_k}^{t_{k+1}}H(s,M_{t_{k+1}})ds\right]\geq  0.
\end{equation}
Summing up the right-hand side of (\ref{wH}) over all $k\in\{ 0,\ldots,n-1\}$ and recalling that $w(T,\cdot)=0$, we get
\[
  w(t_0,m_0)\leq E\left[\sum_{k=0}^{n-1} \Int_{t_k}^{t_{k+1}}H(s,M_{t_{k+1}})ds\right]
\]
Since $(M_s)$ is c\`adl\`ag, letting $\tau\to 0$ yields to 
\[
  w(t_0,m_0)\leq E\left[ \Int_{0}^{T}H(s,M_{s})ds\right]
\]
Now we take  the infimum over $M$ to get the desired result: 
\[ w(t_0,m_0) \leq  \bV^+(t_0,m_0) 
 \]
\QED

\subsection{Characterization of the value of the game as dual solution of a HJI-equation.}

The main drawback of the previous characterization is that it is of nonlocal nature. We now show that it is possible to characterize the value function 
by local inequalities. We use here several ideas of \cite{cls}. 

\begin{Definition} Let  $V:[0,T]\times \Pi_2\to\R$. 
\begin{enumerate}
\item {\em Dual supersolution:} We say that $V$ is a viscosity dual supersolution to 
\be\label{HJdual}
\frac{\partial V}{\partial t}(t,m) + H(t,m)= 0 
\ee
if $V$ is  lower semicontinuous, $V(t,\cdot)$ is convex for any $t$ and if, for any $\bar m\in \Pi_2$ and smooth test function $\varphi:[0,T]\to \R$ such that $t\to V(t,\bar m) -\varphi(t)$ has a local minimum at $\bar t\in [0,T)$ with $\bar m$ extreme point of the graph of $V(\bar t,\cdot)$, one has 
$$
\varphi'(\bar t ) + H(\bar t,\bar m)\leq 0 
$$
\item {\em Dual subsolution:}  We say that $V$ is a viscosity dual subsolution to (\ref{HJdual})
if $V$ is upper semicontinuous, $V(t,\cdot)$ is convex for any $t$ and if, for any $\bar m\in \Pi_2$ and smooth test function $\varphi:[0,T]\to \R$ such that 
$t\to V(t,\bar m) -\varphi(t)$ has a local minimum at $\bar t\in [0,T)$, one has 
$$
 \varphi'(\bar t) + H(\bar t,\bar m)\geq 0 
$$
\item {\em Dual solution:} We say that a continuous map $V:[0,T]\times \Pi_2\to\R$ is a viscosity dual solution of (\ref{HJdual}) if it is a sub and a supersolution of (\ref{HJdual}). 
\end{enumerate}
\end{Definition}

\begin{Proposition}\label{Prop:Charac2} The value function $\bV^+$ is a dual  solution to 
\begin{equation}
\label{dual}
\left\{\begin{array}{l}
\frac{\partial V}{\partial t}(t,m) + H(t,m)= 0 \\
V(T,m)= 0\qquad \forall m\in \Pi_2
\end{array}\right.
\end{equation}
\end{Proposition}

\noindent {\bf Proof : } 
By Proposition \ref{largest}, we know already that $\bV^+$ is a subsolution of (\ref{dual}). 
Let now any $\bar m\in \Pi_2$ and  $\varphi:[0,T]\to \R$ be smooth test function such that 
$ \bV^+(t,\bar m) \geq \varphi(t)$ with an equality at $\bar t$ and $\bar m$ is an extreme point of the graph of $\bV^+(\bar t,\cdot)$. Let $\tilde M$ be the optimal martingale measure as in Corollary \ref{PGD}: 
\begin{equation}
\label{dppm}
\bV^+(\bar t,\bar m)= E\left[\Int_{\bar  t}^{\bar t+h}H(s,\tilde M_s)ds+ \bV^+(\bar t+h, \tilde M_{\bar t+h})\right]\qquad \forall h\in [0,T-\bar t]\;.
\end{equation}
Taking $h=0$ in (\ref{dppm}) gives
$$
\bV^+(\bar t,\bar m)= E\left[\bV^+(\bar t, \tilde M_{\bar t^+})\right]\;, 
$$
so that $\tilde M_{\bar t^+}=\bar m$ because $\bar m$ is an extreme point of $\bV^+(\bar t, \cdot)$. 
Since $\bV^+$ is convex in $m$, we have
\[ E[\bV^+(\bar t+h,\tilde M_{\bar t+h})]\geq \bV^+(\bar t+h,\bar m).\]
Replacing this in the right hand side term of (\ref{dppm}), we obtain
\begin{equation}
\label{fifi}
 \varphi(\bar t)=\bV^+(\bar t,\bar m)\geq  E\left[\int_{\bar t}^{\bar t+h}H(s,\tilde M_s)ds\right]+\bV^+(\bar t+h,\bar m)\geq 
 E\left[\int_{\bar t}^{\bar t+h}H(s,\tilde M_s)ds\right]+\varphi(\bar t+h)
.
\end{equation}
Dividing by $h$ and letting $h\to 0^+$ in the above inequality gives the desired result. \QED

\begin{Proposition} [Comparison] If $V_1$ is a dual supersolution and $V_2$ is a dual subsolution, with $V_1(T\cdot)\leq V_2(T,\cdot)$  in $\Pi_2$, then 
$V_1\leq V_2$ in $[0,T]\times \Pi_2$. 
\end{Proposition}

\begin{Remark} In particular $V$ is the unique dual viscosity solution of equation (\ref{dual}).
\end{Remark}

\noindent {\bf Proof : } Let us fix $\eta, M>0$. We argue by contradiction, assuming that, for $\epsilon>0$ sufficiently small,  
$$
\sup_{t,s,m} V_1(s,m)-V_2(t,m) -\frac{(s-t)^2}{2\epsilon}-\eta t\; >\; 0\;,
$$
where the supremum is taken over the Borel probability measures $m$ with a support in $B(0,M)$. Since this set is compact, there is a maximum point
$(\bar s^\epsilon,\bar t^\epsilon, \bar m^\epsilon_0)$. As $V_1(T\cdot)\leq V_2(T,\cdot)$, we have $\bar s^\epsilon,\bar t^\epsilon<T$ when $\epsilon$ is small enough. Let $C_\epsilon$ be the set of maximum points of the form $(\bar s^\epsilon,\bar t^\epsilon, m)$ and let $\bar m^\epsilon$
be an extreme point of $C_\epsilon$. By Carath\'eodory Theorem, $\bar m^\epsilon$ belongs to $C_\epsilon$. We claim that $\bar m^\epsilon$ is an extreme point of the graph of 
$V_2(\bar t^\epsilon,\cdot)$. Indeed, if $\bar m^\epsilon= \frac{m_1+m_2}{2}$ with $V_2(\bar t^\epsilon,m_1)=V_2(\bar t^\epsilon,m_2)= V_2(\bar t^\epsilon,\bar m^\epsilon)$, we have 
$$
\frac12\sum_{i=1}^2 V_1(\bar s^\epsilon,m_i)-V_2(\bar t^\epsilon,m_i) -\frac{(\bar s^\epsilon-\bar t^\epsilon)^2}{2\epsilon}
\geq V_1(\bar s^\epsilon,\bar m^\epsilon)-V_2(\bar t^\epsilon,\bar m^\epsilon) -\frac{(\bar s^\epsilon-\bar t^\epsilon)^2}{2\epsilon}
$$
by convexity of $V_1(\bar s^\epsilon,\cdot)$. Now note also that the support of $m_1$ and the support of $m_2$ are in $B(0,M)$ because so is the support of $\bar m^\epsilon$. By optimality of $\bar m^\epsilon$ the points $(\bar s^\epsilon,\bar t^\epsilon,m_1)$ and $(\bar s^\epsilon,\bar t^\epsilon,m_2)$ belong to $C_\epsilon$ and therefore $m_1=m_2=\bar m^\epsilon$ since $\bar m^\epsilon$ is an extreme point. 

We now use the definition of viscosity dual solutions: since the  map $s\to V_1(s,\bar m^\epsilon)-\frac{(s-\bar t^\epsilon)^2}{2\epsilon}$ has a maximum
at $\bar s^\epsilon$, we get 
$$
\frac{\bar s^\epsilon-\bar t^\epsilon}{\epsilon}+ H(\bar s^\epsilon, \bar m^\epsilon) \geq 0\;.
$$
Since the  map $t\to V_2(t,\bar m^\epsilon)+\frac{(\bar s^\epsilon-t)^2}{2\epsilon}-\eta t$ has a minimum at $\bar s^\epsilon$ with $\bar m^\epsilon$ an extreme point of the graph of 
$V_2(\bar s^\epsilon,\cdot),$ we get 
$$
\eta+ \frac{\bar s^\epsilon-\bar t^\epsilon}{\epsilon}+ H(\bar s^\epsilon, \bar m^\epsilon) \leq 0
$$
So
$$
\eta+ H(\bar s^\epsilon, \bar m^\epsilon)-H(\bar s^\epsilon, \bar m^\epsilon)\leq 0\;.
$$
As $\epsilon\to 0$, we have $\bar s^\epsilon, \bar t^\epsilon\to \bar t$ and $\bar m^\epsilon\to \bar m$ (up to some subsequence), so that $\eta<0$. This contradicts the definition of $\eta$. 

In conclusion,  for any $\epsilon>0$ and $\eta>0$ and for any $(s,m)\in [0,T]\times \Pi_2$ such that the support of $m$ is in $B(0,M)$, we have
$$
V_1(s,m)-V_2(t,m) -\frac{(s-t)^2}{2\epsilon}-\eta t\; \leq\; 0\;,
$$
Taking $s=t$, and letting $\eta\to 0^+$ and $M\to +\infty$, we obtain that inequality
$V_1(t,m)-V_2(t,m) \; \leq\; 0$ holds for any $(t,m)$ such that $m$ has a compact support. We complete the proof by density. 
\QED

\bigskip

\noindent{\bf Acknowledgments}\\
This work has been supported by the Commission of the European Communities under the 7-th Framework Programme Marie Curie Initial Training Networks Project ``Deterministic and Stochastic Controlled Systems and Applications'' FP7-PEOPLE-2007-1-1-ITN, no. 213841-2 and the French National Research Agency ANR-10-BLAN 0112.

\end{document}